%%%   3 May 1995
%%%%%%%%%%%%%%%%%%%%%%%%%%%%%%%%%%%%%%%%%%%%%%%%%%%%%%%%%%%%%%%%%%%%%%%%%%  
%%%
%%%     Universal Non-Completely-Continuous Operators 
%%% 
%%%
%%%          Maria Girardi and  William B. Johnson
%%%
%%%
%%%                       AMS-TeX 2.1 
%%%
%%%%%%%%%%%%%%%%%%%%%%%%%%%%%%%%%%%%%%%%%%%%%%%%%%%%%%%%%%%%%%%%%%%%%%%%%%

%**%Before printing, find \voffset command

\documentstyle{amsppt}
\tolerance=30000
\NoBlackBoxes
\magnification=\magstep1 
%\NoRunningHeads
\vsize=7.4in     %**%to get page nos.
%\voffset=-.4in  %**%needed for getting page no. on page 1 on some machines
%**% next line may be needed on some machines
%\input mssymb %%%This and \def\frak{\Cal} is what worked best

%%%%%%%%%%%%%%%%%%%%%%%%%%% MACROS   %%%%%%%%%%%%%%%%%%%%%%%%%%%%%%%

%**% If you don't have \frak, use next line:
%\def\frak{\Cal}

\def \qed {\vrule height6pt  width6pt depth0pt} %%%End of proof symbol
\def\X{\frak X} 
\def\Y{\Cal Y} 
\def\Z{\Cal Z} 
\def\D{\Cal D}
\def\C{\frak C}
\def\fS{\frak S}
 
\def\td{\tilde}
\def\hf{\hfill}
\def\d{\delta}
\def\a{\alpha} 
\def\e{\epsilon}
\def\Sg{\Sigma}
\def\S+{\Sigma^+} 
\def\n{\Vert} 
\def\av{\mid}

\def \an{\{ a_n \}_{n=1}^{\infty}}

%%%%%%%%%%%%%%%%%%%%%%%%%%%%%%%%%%%%%%%%%%%%%%%%%%%%%%%%%%%%%%%%%%%%%%%%%%  

\topmatter  
\title   Universal Non-Completely-Continuous Operators 
\endtitle 
\author  Maria Girardi 
\footnote "*"{Supported in part by NSF grants DMS-9306460 and DMS-9003550. 
           \hfill\phantom{P}}
\footnote"{$^\ddagger$}"{Participant, NSF
Workshop in Linear Analysis
\& Probability, Texas A\&M University. \hfill\phantom{M}}
\quad and \quad William B. Johnson$^*$
\endauthor
\date 
    3 May 1995  
\enddate
\subjclass 
47B99, 47A68, 47D50, 47B38, 
46B20, 46B28,  46B07, 46B22
\endsubjclass 
 
\abstract {
A  bounded  linear operator  between Banach spaces  is called {\it
completely continuous} if it carries weakly   convergent sequences into
norm convergent  sequences.   Isolated is a universal operator  for the
class of  non-completely-continuous operators from $L_1$ into an 
arbitrary Banach space, namely,  the   operator  from $L_1$ into
$\ell_\infty$  defined by   
$$  T_0 (f) =\left( \int r_n f \, d\mu \right)_{n\ge 0} \  , 
$$   where $r_n$ is the $n^{\text{th}}$ Rademacher function.   It is
also shown that there does not exist  a  universal operator for the
class of  non-completely-continuous operators between two  arbitrary
Banach space.  The proof  uses the factorization theorem for weakly
compact  operators and a Tsirelson-like space.}
\endabstract      
\endtopmatter
\address  Maria Girardi
\newline University of South Carolina,  Department of Mathematics, 
Columbia, SC  29208 \endaddress     \email
girardi\@math.scarolina.edu  \endemail

\address W.~B.~Johnson
\newline Department of Mathematics, Texas A\&M University, College
Station, TX 77843 \endaddress
\email johnson\@math.tamu.edu\endemail

\headline={\ifodd\pageno\rightheadline\else\leftheadline\fi}
\def\rightheadline{\eightrm\hfil UNIVERSAL OPERATORS \hfil\folio}
\def\leftheadline{\eightrm\hfil GIRARDI, JOHNSON \hfil\folio}

\document 
\baselineskip 15 pt

Suppose that  $\C$ is a class of (always bounded, linear, between 
Banach spaces) operators so that an
operator
$S$ is in {$\C$} whenever the domain of $S$ is the domain of some operator in
{$\C$} and there exist operators $A$, $B$ so that $BSA$ is in {$\C$}; the natural
examples of such classes are all the operators that do not belong to a given
operator ideal.  A subset  {$\fS$} of such a class
$\C$ is said to be {\it universal}\/ for $\C$ provided for each $U$ in $\C$, some
member of
$\fS$ factors through $U$; that is, there exist operators $A$ and $B$ so that
$BUA$ is in $\fS$.  In case $\fS$ is singleton; say, $\fS=\{S\}$; we say that $S$
is universal for $\C$.  

In order to study a class $\C$ of operators, it is natural
to try to find a universal subclass of $\C$ consisting of specific, simple
operators.  For certain classes, such a subclass is known to exist.  For
example, Lindenstrauss and Pe\l czy\'nski, who introduced the concept of
universal operator, proved  [LP] that the ``summing operator" from $\ell_1$ to
$\ell_\infty$, defined by \ $\an\mapsto \{\sum_{k=1}^n a_k\}_{n=1}^\infty$, is
universal for the class of non-weakly-compact operators; while in [J] it was
pointed out that the formal identity from $\ell_1$ to
$\ell_\infty$ is universal for the class of non-compact operators.

An operator  between Banach spaces 
is called {\it completely continuous} if it carries weakly  
convergent sequences into norm convergent  sequences.    
The operator   from $L_1$ into $\ell_\infty$ given by  
$$ 
T_0 (f) =\left\{ \int r_n f \, d\mu \right\}_{n= 0}^\infty \  , 
$$  
where $r_n$ is the  $n^{th}$ Rademacher function, 
is not completely continuous.   
We prove in Corollary 4 that  $T_0$ is universal for the class of
non-completely-continuous operators from an $L_1$-space; 
however,  in Theorem 5  we show that
there does not exist a universal  non-completely-continuous operator.  

Throughout this paper, $\X$ denotes an arbitrary Banach space, 
$\X^*$ the dual space of~$\X$,   
and  $S(\X)$ the unit sphere of $\X$.   
The triple $(\Omega, \Sigma, \mu)$~refers 
to the Lebesgue measure space on  $[0,1]$, 
$\S+$ to the sets in $\Sigma$ with positive measure, 
and $L_1$ to $L_1(\Omega, \Sigma, \mu)$.
All notation and terminology, not otherwise explained, are as in [DU] 
or [LT]. 

To crystalize the ideas in Theorem~1,  we introduce some
terminology.    A system 
$\Cal A = \{ A^n_k \in\Sg \: n=0,1,2,\ldots \text{ and } k=1,\ldots,2^n \}$  
is a {\it dyadic splitting} of $A^0_1\in\S+$ if  
each $A^n_k$  is partitioned into the   
two sets $A^{n+1}_{2k-1}$ and $A^{n+1}_{2k}$ 
of equal measure for each admissible $n$ and $k$ .  
Thus  the collection $\pi_n = \{ A^n_k \: k=1, \ldots, 2^n \}$ of sets 
along the $n^{\text{\, th}}$-level  partition $A^0_1$ with 
$\pi_{n+1}$ refining $\pi_n$ and $\mu(A^n_k) = 2^{-n}\mu(A^0_1)$. 
To a dyadic splitting corresponds 
a (normalized) Haar system $\{ h_j \}_{j\ge 1}$  along with its   
natural blocking $\{ H_n \}_{n \ge 0}$  
where 
$$
   h_1 = \tfrac{1}{\mu(A^0_1)} \, 1_{A^0_1}
   \text{\qquad and \qquad} 
  h_{2^n+k} = \tfrac{2^n}{\mu(A^0_1)} \, 
               (1_{A^{n+1}_{2k-1}} - 1_{A^{n+1}_{2k}} ) 
$$ 
for $n=0,1,2,\ldots $,  $k=1,\ldots,2^n $,    
and  $H_n = \{ h_j \:  2^{n-1} < j \le 2^n \}$.   
The usual Haar system $\{ \tilde h_j \}$ 
corresponds to the usual dyadic splitting 
$\left\{ \left[ \frac{k-1}{2^n}, \frac{k}{2^n} \right) \right\}_{n,k}$.  
Let $L_1 (\Cal A)$ be the closed subspace of $L_1$ with basis 
  $\{h_j \}_{j\ge 1}$.

A set $N$ in the  unit sphere of the dual of a 
Banach space $\X$ is said to norm  a subspace $\X_0$ 
within $\tau >1$ if for each $x\in\X_0$  there is $x^*\in N$ 
such that $\n x \n \le \tau x^* (x)$. 
It is well known and easy to see that  a   
sequence $\{\X_j\}_{j \geqslant 1}$ of subspaces of $\X$ forms a 
finite dimensional decomposition 
with constant at most $\tau$ provided that  for each $n\in\Bbb N$ 
the space  generated by $\{ \X_1, \ldots, \X_n \}$ 
can be normed by a  set from $S(\X_{n+1}^\perp)$ within $\tau_n > 1$ 
where $\Pi\tau_n\le\tau$.    

To help demystify   
Theorem~1, 
we examine more closely  the 
operator  $T_0\: L_1 \to \ell_\infty$ given above. 
This operator does more than 
just map the Rademacher functions   
$\{ r_n \}$ to the standard unit vectors  $\{ e_n \}$ in  
$\ell_\infty$ (which suffices to guarantee that it is   
not completely continuous).  Let $ x^*_n$ 
be the $n^{\text{th}}$ unit vector of $\ell_1$, 
viewed as an element in the dual of~$\ell_\infty $.  
For the usual dyadic splitting of the unit interval, 
$r_n$ is just the sum of the Haar functions in $H_n$, 
properly normalized.  
Thus $1 = \n T_0 r_n \n = x^*_n (T_0 r_n )$ follows 
from the stronger condition that 
$$ 
x^*_n (T_0 h) =  \d_{n,m}  \text{\qquad for each\quad} h \in H_m  \ .
$$  
Note that $ T_0^*x^*_n  $ is just $r_n$, which as  
a sequence in $L^*_1$ is weak*-null and  
equivalent to the unit vector basis of $\ell_1$. 
Since $T_0$ maps each element in $H_n$ to $e_n$, 
the collection $\{ \text{sp } T_0 H_n \}$   forms 
a finite dimensional decomposition. 
Theorem~1 states that each   
non-completely-continuous operator $T$ on $L_1$ 
behaves like the 
operator  $T_0$ in the sense that there is 
some dyadic splitting of some subset of $[0,1]$  
so that the corresponding Haar system with  
 $T$ enjoys the above properties of the usual Haar system 
  with $T_0$.

\proclaim{Theorem~1} 
Let $Y$  be a subset of  $S(\X^*)$  that 
norms $\X$ within some fixed constant greater than one 
and  let $\Y$ be a subspace   of $\X^*$ that contains $Y$.    
If the  operator  $T\: L_1\to\X$ is not completely continuous 
and $\{ \tau_n \}_{n\ge 0}$  is a sequence of numbers larger than 1, 
then there exist 
\roster 
\item"{(A)}" 
a dyadic splitting $\Cal A = \{ A^n_k \}$ 
\item"{(B)}"
a sequence $\{ x^*_n\}_{n\ge0}$ in $S(\X^*) \cap \Y$ 
\item"{(C)}"
a finite set $\{ z^*_{n,i} \}_{i=1}^{p_n}$ in  $S(\X^*)$ for  each $n\ge 0$ 
\endroster 
such that for the Haar system  $\{ h_j \}_{j\ge 1}$  and  the 
blocking $\{ H_n \}_{n \ge 0}$   
corresponding to $\Cal A$, for some $\delta>0$, and each $n,m \ge 0$, 
\roster 
\item $x^*_n (T h) = \d \cdot \d_{n,m} $ \quad for each $h \in H_m $  
\item $\{T^* x^*_n \} $ \quad   is weak*-null in $L_\infty$ 
\item  $\{ T^* x^*_n \}$  \quad is   
    equivalent to the unit vector basis of $\ell_1$ 
\item $\{ z^*_{n,i} \}_{i=1}^{p_n}$  \quad norms  
      $\text{sp} (  \cup_{j=0}^n T H_j )$ 
       within $\tau_n$ 
\item  $T H_{n+1} \subset  \  ^\perp \{z^*_{n,i}  \}_{i=1}^{p_n} $  . 
\endroster  
Note that condition {\rm  (3)} implies that    
$\{x^*_n \}$ is also  
equivalent to the standard unit vector basis of $\ell_1$.  
If $\Pi \tau_n $ is finite, 
then the last two conditions     
guarantee that   $\{ \text{sp } T H_n \}_{n \ge 0}$  
forms a finite dimensional decomposition  with constant at most $\Pi \tau_n $.
\endproclaim

The proof    uses the following   two 
standard  lemmas.  
 
\proclaim{Lemma~2}  
Let $E = \text{sp } \{ x_i \}_{i=0}^m$ be a finite dimensional subspace of a
Banach space $\X$ and let $\Y$ be a total subspace of $\X^*$. 
For each $\e >0$ there exists
$\eta >0$    such that if   $y^* \in \X^*$ satisfies 
$\left| y^* (x_i) \right| < \eta$  
for each 
$ 1 \le i \le m$, then 
there exists  $x^* \in  E^\perp$  
of norm $0$ or $\n y^* \n$ 
such that $\n x^* - y^* \n < \e$.  
Furthermore, if $y^*$ is in $ \Y$         
then $x^*$ can be taken to be in $\Y$.   
\endproclaim 

\demo{Proof of Lemma~2}  
Assume, without loss of generality, 
$ \{ x_i \}_{i=0}^m$ is linearly independent. 
Consider the isomorphism 
$l \: E \to \ell_1^m$ that takes $x_i$ to the $i^{\text{\, th}}$  
unit basis vector of $\ell_1^m$ 
and  let  $P$ be a projection from $\X$ onto $E$ that is
$w(\Y)$-continuous, so that $P^*E^*$ is a subspace of $\Y$.  Such a
projection exists because $\Y$ is total. Then 
$\tilde x^*
\equiv y^*
\cdot
\left( I_{\X} - P
\right)$   is in   $E^\perp$ .  
It is easy to check that for 
$\eta = \frac{\e}{3~\n l \n~ \n P \n}$,  
if $\left| y^* (x_i) \right| < \eta$ for each $i$, 
then  $\n \tilde x^* - y^* \n \le \frac{\e}{3}$. 
If  $\n \tilde x^* \n = 0$, then let $x^* = \tilde x^*$. 
Otherwise, let  
$x^* =\left( \n y^* \n ~/~ \n \tilde x^* \n\right)~\tilde x^* $.   
Then $\n x^* - y^* \n  \le  2 \n \tilde x^* - y^* \n $. 
Thus $x^*$ does what it is to do.   
\hfill\qed
\enddemo   
  
Recall that the extreme points of $B(L_\infty)$ are just the 
$\pm 1$-valued measurable functions. 

\proclaim{Lemma~3}  
If $\{ f_i \}_{i=0}^n$ is a finite subset of $L_1$, 
$\{ \a_i \}_{i=0}^n$ are scalars,  and 
$$
  S = \left\{ g \in B(L_\infty) \:  \int f_i g \, d\mu = \a_i 
      \text{ \ for each \  }  0 \le i \le n  \right\} \ ,  
$$  
then $\text{ext }   S = S \cap   \text{ext }  B(L_\infty)$, 
where  ext denotes the extreme points of a set.  
Also,  if $S$ is non-empty then so is $ \text{ext } S$. 
\endproclaim 
\flushpar 
Specifically, we use the following version 
of this extreme point argument lemma. 
\proclaim{Lemma~3$^\prime$}
If  $F = \{f_1, \ldots, f_n \}$  
and  there exists $g$ in $ B(L_\infty)\cap F^\perp$  
such that  
$\int f_0 g \, d\mu \ge  \alpha_0 > 0$, 
then there exists a  $\pm 1$-valued  function 
$u$ in $ B(L_\infty)\cap F^\perp$  
such that 
$\int f_0  u \, d\mu = \alpha_0 $. 
\endproclaim 
 
\demo{Proof of Lemma~3}  
Consider, if there is one, a  function $g$ in  $S$ 
for which there exists a subset $A$ of positive measure and $\e >0$ 
such that $  -1 + \e < g 1_A < 1-\e$.  
Since the set $ \{ f \in L_\infty \: \mid f \mid~\le  1_A \} 
   \cap  \{ f_i \}_{i=0}^n\,^\perp$ \,  
is infinite dimensional, it contains a non-zero element $h$ of  
norm less than $\e$.  But then $g \pm h \in S $ and 
 so $g$ is not an extreme point of $S$.  
Thus   $ \text{ext } S = S \cap  \text{ext }  B(L_\infty)$.   

Since $S$ convex and weak*-compact in $L_\infty$, 
if $S$ is non-empty then so is $\text{ext } S$. 

As for the last claim of the lemma, 
just note that if 
$g\in B(L_\infty)\cap F^\perp$  satisfies  
$\int f_0 g \, d\mu \equiv  \beta \ge \alpha_0 > 0$,  
then $  \frac{\alpha_0}{\beta} ~ g$ is in  the set $S$ 
where   $\alpha_i = 0$ for $i >0$.  
By the first part of the lemma, 
any extreme point $u$ of $S$  will do. 
\hfill\qed 
\enddemo 

Although the proof of Theorem~1 is somewhat technical, 
the overall idea is simple.  
Since $T$ is not completely continuous,  
we start by finding a weakly convergent sequence 
$\{ g_n \}$ in $L_1$ and  norm one functionals $y^*_n$ 
such that  $\delta_0 \le y^*_n  \left(T \, g_n\right)$. 
 Each $x^*_n$ will be a small perturbation of some $y^*_{j_n}$.
Conditions (2) and (3) can 
 be arranged by standard arguments.

Now the proof gets technical.  
We begin by finding a subset $A^0_1$ where the 
$L_\infty$ function 
$(T^*y_n^*) g_n$, 
which in the motivating  example of $T_0$ is the function $r_n r_n$, 
is large in some sense.  We then proceed by 
induction on the level $n$.  
Given a finite dyadic splitting 
up to $n^{\, \text{th}}$-level 
provides  the  subsets $\{ H_m \}_{m=0}^n$ of corresponding 
Haar functions.   We need to split each $A^n_k$ into 2 sets 
$A^{n+1}_{2k-1}$ and  $A^{n+1}_{2k}$ (thus finding $h_{2^n+k}$) 
and find the desired functionals so that all works.  
It is easy to find the functionals  to satisfy condition~(4).  
In the search for $x^*_{n+1}$, 
apply Lemma~2 to the set $E$ given in $(\dagger)$ 
 so that 
we need only to almost (within some $\eta$) 
 satisfy (1-i$^\prime$)  for some  
$y^*_{j}$;  
for then we can perturb $y^*_{j}$ to find $x^*_{n+1}$ 
that satisfies (1-i$^\prime$) exactly.  
Next, for each $A^n_k$, apply Lemma~3$^\prime$ 
with $F$ as given in $(\ddagger)$ and $f_0 = T^* y_j^* 1_{A^n_k}$ 
and $g$ being a small perturbation of $g_j 1_{A^n_k}$.
All is set up so that
such a perturbation exists   for 
a $j$ (dependent on $n$ but independent of $k$) 
 sufficiently  large enough.  
Now Lemma~3$^\prime$ gives that desired $\pm 1$-valued  Haar-like function 
that yields the desired splitting of the 
 $(n+1)^{\, \text{th}}$-level.   
The sets $F^n_k$ are chosen exactly so that 
conditions (1-ii$^\prime$), (1-iii$^\prime$), and 
(5$^\prime$) hold.

\demo{Proof of Theorem~1} 
Let $T\: L_1\to\X$ be a norm one  operator  that is not completely continuous. 
Then there is a sequence $\{ g_n \}$ in $L_1$ 
and  a sequence $\{ y^*_n \}$ in $S(\X^*) \cap \Y$ satisfying: 
\roster
\item"{(a)}" $\n g_n \n_{L_\infty} \ \le 1$ 
\item"{(b)}" $g_n$  \  is weakly null in  $L_1$ 
\item"{(c)}" $\delta_0 \le y^*_n  \left(T \, g_n\right)$  \ for some $\d_0
> 0$ .
\endroster  
Using (a), (b), and (c) along with Rosenthal's \ $\ell_1$ \ theorem
[cf\.~LT,~Prop\.~2.e.5],  by passing 
to a further subsequence, we also have  that  
\roster
\item"{(d)}" $\{ T^* y_n^* \}$   \  is 
  equivalent to the standard unit vectors basis of $\ell_1$.
\endroster  
Since $B(L_\infty)$ is    weak* sequentially-compact 
in $L_\infty$,   
by passing to a subsequence and considering 
differences we may assume  that 
\roster 
\item"{(e)}"  $T^* y^*_n$  \  is weak*-null in $L_\infty$, 
\endroster   
where (d) allows normalization of the new $y^*_n$'s so as to keep  them 
in $S(X^*)$ and, used with care, 
(b) ensures that (c) still holds for  some (new) positive $\d_0$.  
But  $\{\left( T^* y^*_n \right)\cdot  g_n\}$ is also in $B(L_\infty)$ and
so,  by passing to yet another subsequence, we have that 
\roster
\item"{(f)}"  $\{\left( T^* y^*_n \right)\cdot g_n\} \to h $ \ weak* in
$L_\infty$ 
\endroster 
for some $h \in L_\infty$.

Since  $\int h \, d\mu \ge \d_0$,   
the set $A \equiv  [h\ge\d_0]$ has positive measure. 
We may assume, by replacing $y^*_n$ by $-y^*_n$ 
and $g_n$ by $-g_n$ when needed, that 
$\n T^* y^*_n \mid_A \n_{L_\infty} = \text{ess sup } T^* y^*_n \mid_A$. 
So from (a) and (f) it follows that 
$\delta_0 \leqslant \liminf  \text{ess sup } T^* y^*_n \mid_A$ 
while from (e) it follows that 
$\limsup \mu [   T^* y^*_n \mid_A \geqslant \delta_0 - \eta ] < \mu(A)$ 
for each $0< \eta < \delta_0$.  
Thus, since the closure of the set 
$$  
    \left\{   \frac{\int_E f \, d\mu}{\mu(E)} 
    \: E \subset A , E\in\Sigma^+ \right\} 
$$     
is the interval 
$[ \text{ess inf } f , \text{ess sup } f]$,  
there is a subset $A^0_1$ of $A$  with positive measure  
and $j_0$ such that 
$y^*_{j_0} T(1_{A^0_1}) = \d \mu(A^0_1)$ 
for some positive $\d$  less than $\d_0$, say  $\d  \equiv \d_0 - 3\e $.
Put $x_0^* = y^*_{j_0}$ 
and $H_0 \equiv \{ h_1 \} = \left\{ 1_{A^0_1} / \mu(A^0_1) \right\}$.

We shall construct, by induction  on the level $n$, 
a dyadic splitting of $A^0_1$ along with 
the desired functionals.  
Towards this, take  
a  decreasing sequence $\{ \e_n \}_{n\ge 0}$ of positive numbers  
such that  $\e_0 < \e$   and $\sum \e_n < \frac{\d_0}{2K}$  where 
$K$ is the basis constant of $\{ T^* y_n^* \}$. 
The sequence $\{ x^*_n \}$ will be chosen 
such that $\n x^*_n - y^*_{j_n} \n \le  \e_{n}$ 
for some increasing sequence $\{ j_n \}_n$ of integers,  
which will ensure   
conditions~(2)~and~(3).  
Note that condition (1) is equivalent to 
the following 3 conditions  holding  
\roster
\item"{(1-i)}"   
$x^*_{n} (T h) = 0 $ \quad for $h\in H_m$ and  $0 \le m <  n$  
\item"{(1-ii)}"  
$x^*_m (T h) = 0 $ \quad for $h\in H_{n}$ and $0 \le m  <  n$ 
\item"{(1-iii)}"  
$x^*_{n}(T h) = \d $ \quad for $h \in H_{n}$  
\endroster  
for each $n$.  
Clearly these  three conditions hold for $n = 0$.    
Fix $n\ge0$.

Suppose that we are given a finite dyadic splitting  
$\{ A^m_k \: m=0,\ldots, n \text{ and } k=1,\ldots,2^m \}$  of $A^0_1$
up to $n^{\, \text{th}}$-level, 
which    gives the  subsets $\{ H_m \}_{m=0}^n$ of corresponding 
Haar functions.  
Thus  we can find a finite set  $\{ z^*_{n,i} \}_{i=1}^{p_n}$ in  $S(\X^*)$ 
such that   
$\{ z^*_{n,i} \}_{i=1}^{p_n}$ norms  
      $\text{sp} (  \cup_{j=0}^n T H_j )$ 
       within $\tau_n$. 
Suppose that we are also given  $\{ x^*_m \}_{m=0}^n$   in
$\Y\cap S(X^*)$ such that
the three subconditions of $(1)$ hold
and 
if $k=1,2,\dots,n$ then $ \n x^*_k - y^*_{j_k} \n \le \e_k$  for
some  $j_k$.  

We shall find $x^*_{n+1}$ along with 
$j_{n+1} > j_n$ such $\n x^*_{n+1} - y^*_{j_{n+1}} \n \le \e_{n+1}$  
and we shall partition, for each  $ 1\le k \le 2^n$, 
the set $A^n_k$ into 2 sets 
$A^{n+1}_{2k-1}$ and  $A^{n+1}_{2k}$ of equal measure 
(thus finding $h_{2^n+k}$ and so finding   the corresponding set 
$\{ H_{n+1} \}$) such that    
\roster
\item"{(1-i$^\prime$)}"  
$x^*_{n+1} (T h) = 0 $ \quad for $h\in H_m$ and $0 \le m  < n+1$
\item"{(1-ii$^\prime$)}"  
$x^*_m (T h) = 0 $  \quad for $h\in H_{n+1}$ and $0 \le m  < n+1$ 
\item"{(1-iii$^\prime$)}"  
$x^*_{n+1} (T h) = \d $  \quad for $h \in H_{n+1}$  
\item"{(5$^\prime$)}"  
   $T H_{n+1} \subset \  ^\perp \{z^*_{n,i}  \}_{i=1}^{p_n} $ .   
\endroster  

Towards this, apply  
Lemma~2 to 
$$
       E \equiv \{ T h \: h \in  H_m \ , \ 0 \le m \le n   \}
       \tag{$\dagger$}
$$ 
and $\e_{n+1}$ to find the corresponding $\eta_{n+1}$. 
Let 
$$
F^n_k = \{ 1_{A^n_k}\} \cup  
        \{T^*x_m^*  1_{A^n_k}\}_{m=0}^n \cup   
        \{T^*z_{n,i}^*  1_{A^n_k}\}_{i=1}^{p_n}   
        \subset L_1   \tag{$\ddagger$}
$$
and $F_n = \text{sp } \left[ ~\cup_{k=1}^{2^n} F^n_k~\right] $.  
  
Pick $j \equiv j_{n+1} > j_n$ so large that for  
$k=1,\ldots, 2^n$ 
\roster	 
\item"{(g)}" $\left| \left(T^*y_j^* \right) h \right| < \eta_{n+1} $
   \quad for all $h\in \cup_{m=0 }^n H_m$ 
\item"{(h)}"$\left|  \int_\Omega g_j f \, d\mu \right| 
        \le  {\frac{\e}{3}}   \n f \n $ \quad  for all $f$ in $F_n$     
\item"{(i)}"$\int_{A^n_k} T^*y^*_j \cdot  g_j \, d\mu  
        \ge \left(\d_0 - \e \right)  \mu(A^n_k)$ . 
\endroster   
Condition (g) follows from (e),  
condition (h)  follows from (b) and the 
fact that $F_n$ is finite dimensional, 
condition (i) follows from  (f) and  
the definition of $A$.   
 
By Lemma~2 and   (g), 
there  is  $x^*_{n+1} \in  S(\X^*) \cap \Y$ such that 
$\n x^*_{n+1} - y^*_{j_{n+1}} \n$ is at most $ \e_{n+1}$  
and  $x^*_{n+1} Th = 0$ for each $h \in \cup_{m= 0}^n H_m$.  
Thus (1-i$^\prime$) holds.  

Condition (h) gives that 
the  $L_\infty$-distance
from $g_j$ to 
$F_n^\perp  \equiv 
\{ g\in L_\infty \: \int_\Omega f g \, d\mu = 0 \text{ for each } f\in F_n \}$ 
is at most $\tfrac{\e}{3} $.   
So there is $\tilde g_j \in  {F_n^\perp} \cap B(L_\infty)$ such that 
$\n \tilde g_j - g_j \n_{L_\infty}$ is less than $\e $.  
Clearly $\tilde g_j 1_{A^n_k} \in  F^n_k\,^\perp \cap B(L_\infty)$ 
for each admissible $k$.  
By condition~(i),   for each admissible $k$, 
$$ 
\int_\Omega  \left(T^* x^*_{n+1} \right) \cdot
     \left( \tilde g_j  1_{A^n_k} \right)\, d\mu \ge \d  \mu(A^n_k) 
$$ 
and so, by Lemma~3,   
there exists a function  \  $u^n_k \in B(L_\infty) \cap F^n_k\,^\perp$  \    
such that 
$$ 
\int_\Omega \left(T^* x^*_{n+1} \right) \cdot
      \left(  u^n_k  \right)\, d\mu = \d  \mu(A^n_k)  
\tag{$\ast$} 
$$ 
and $u^n_k$ is of the form $1_{A^{n+1}_{2k-1}} - 1_{A^{n+1}_{2k}}$ 
for 2 disjoint sets $A^{n+1}_{2k-1}$ and $A^{n+1}_{2k}$ 
whose union is $A^n_k$.  
Furthermore, $A^{n+1}_{2k-1}$ and $A^{n+1}_{2k}$ are of 
equal measure since $1_{A^n_k} \in F^n_k$. 
Since $u^n_k \in  F^n_k\,^\perp$, 
conditions (1-ii$^\prime$) and (5$^\prime$) hold. 
Condition (1-iii$^\prime$) is just ($\ast$).  
\hfill\qed 
\enddemo

Theorem~1 contains much information.        
For example, the next corollary crystallizes the role of the 
previously mentioned operator $T_0$.   

\proclaim{Corollary~4}    
If the operator $T \: L_1 \to \X$ is not completely continuous, 
then there exist an isometry $A$ and an operator  
 $B$ such that the  following diagram commutes. 
$$ 
\vbox{
\settabs\+\indent
&\hskip 1 true in &\hskip 2 true in &\hskip 1 true in \cr 
\+&\hf$L_1$\hf
  &\hf$\longrightarrow^{^{\hskip-1.25em T}}\hskip1em$\hf&\hf$\X$\hf&\cr
\+\cr
\+&\hf$A\uparrow$\hf&  &\hf$\downarrow B$\hf&\cr
\+\cr 
\+&\hf$L_1$\hf
  &\hf$\longrightarrow^{^{\hskip-1.25emT_0}}\hskip1em$\hf
  &\hf$\ell_\infty$\hf\cr}
$$   
Furthermore, if $\X$ is separable, then 
$T_0$ and $B$  may be viewed as  operators into $c_0$. 
\endproclaim  

\demo{Proof of Corollary~4}   
Let $j_1$ be the natural   injection of 
$L_1(\Cal A)$ into $L_1$,  
let $\X_0$ be the norm closure of $T\left( j_1~ L_1\left(\Cal A\right)\right)$, 
and let  $\td x^*_n$ be the restriction of $x^*_n$ to $\X_0$.  

Since $\{ T^* x^*_n \}$ is weak*-null in $L_\infty$, 
$ \td x^*_n $  is weak*-null in $\X^*_0$. 
Thus the mapping  $U \: \ell_1 \to \X_0^*$  that take the $n^{\text{\, th}}$ 
 unit basis vector of $\ell_1$ to $\td x_n^*$ 
is weak* to weak* continuous and so 
$U$ is the adjoint of the operator 
$S \: \X_0 \to c_0$ where 
$S(x) =  \left( \td x^*_n \left(x\right)\right)_{n \ge 0}$. 

%\vfill\eject  
Consider the (commutative) diagram: 
$$
\vbox{
\settabs\+\indent
&\hskip 1 true in &\hskip 2 true in &\hskip 1 true in 
    &\hskip 1 true in &\hskip 1 true in\cr 
\+&\hf$L_1$\hf
  &\hf$\longrightarrow^{^{\hskip-1.25em T}}\hskip1em$\hf&\hf$\X$\hf&\cr
\+\cr
\+&\hf$j_1\uparrow$\hf&  &\hf$\uparrow j_2$\hf&\cr
\+\cr 
\+&\hf$L_1(\Cal A)$\hf
  &\hf$\longrightarrow^{^{\hskip-1.25em T_{\Cal A}}}\hskip1em$\hf
  &\hf$\X_0$\hf&\cr   
\+\cr   
\+&\hf$R \uparrow$\hf& &\hf$\downarrow S$\hf&\cr
\+\cr
\+&\hf$L_1$\hf&\hf$\longrightarrow$\hf&\hf$c_0$\hf  
   &\hf$\longrightarrow^{^{\hskip-1.25em j_3}}\hskip1em$\hf
   &\hf$\ell_\infty$\hf\cr}
$$  
where  $R \: L_1 \to L_1(\Cal A)$  is the  natural  
isometry that takes a usual Haar function $\widetilde h_j$ in $L_1$ 
to the corresponding associated Haar function  $ h_j$  in $L_1 (\Cal A)$,  
the maps $j_i$ are the natural injections, and 
$T_{\Cal A}$ is such that the upper square commutes.  
  
For an arbitrary space $\X$,  
since $\ell_\infty$  is injective, 
the operator $j_3 S $ extends to an operator 
$\tilde S \: \X \to \ell_\infty$.  
For a separable space $\X$, 
since  $c_0$  is separably injective,   
this extension $\tilde S$ may be view as 
taking values in $c_0$. 
 
Let $A = j_1 R$ and $B = \frac{1}{\delta}\tilde S$. Then   
$B T A ( \widetilde h_{j} )  =  
   \frac{1}{\delta}  
 \left( \tilde x^*_n \left( T h_j \right) \right)_ {n \ge 0}  $.  
Property~1 of Theorem~1 gives that $BTA = T_0$.  
\hfill\qed  
\enddemo

 Corollary~4  says that,   viewed as an operator 
into $\ell_\infty$ (respectively, into $c_0$), 
$T_0$  is universal for the class of  
non-completely-continuous  
operators from $L_1$ into 
an arbitrary (respectively, separable) 
Banach space.

\proclaim{Theorem~5} 
There does not exist  a  universal  operator for the class of 
non-completely-continuous operator.     
\endproclaim  
\flushpar   
The proof of the nonexistence of such an operator  uses 
the existence of a factorization through a reflexive 
space for a  weakly compact operator.  

\demo{Proof}   
Suppose that  there  did exist   
a universal non-completely-continuous operator, say 
$T_1\: \X \to \Z$ where $\X$ and $\Z$ are   Banach spaces.   
Then there is a sequence $\{ x_n \}$ in $\X$ of norm one elements 
that converge weakly to zero  but  
whose images  $\{ T_1 x_n \}$ are  uniformly bounded away from zero. 
Furthermore, by passing to a subsequence, we also have 
that $\{ T_1 x_n \}$ is a basic sequence in $\Z$.    

The first step of the proof uses  $T_1$ to construct 
a   ``nice''  universal non-completely continuous operator.  
By Corollary~7 in [DFJP], there exists 
a reflexive space $\Y$ with  a normalized unconditional basis  
$\{ y_n \}$ such that the map $S \: \Y \to \X$  
that sends $y_n$ to $x_n$ is continuous.   
Consider the map $U\: \Z \to \ell_\infty$ that sends 
$z$ to $( z_n^* (z) )$   
where $\{ z_n^*\}$  is a bounded sequence in  
 $\Z^*$ such that $\{T_1 x_n, z_n^* \}$    
is  a biorthogonal system.   
The map $ I_{\Y} \equiv UT_1S$ sends $y_n$ to the $n^{\text{\, th}}$ 
unit vector of $\ell_\infty$.  
The reflexivity of $\Y$ guarantees that  $I_{\Y }$ is not 
completely continuous.  
Since $I_{\Y }$ factors through the universal 
operator $T_1$,  the operator $I_{\Y }$ must also be universal.  
We now work with this ``nice'' operator $I_{\Y }$. 

For any linearly independent finite set   $\{ x_k \}_{k=1}^n$,     
let $\D \{ x_k  \}_{k=1}^n$  be the norm of the operator 
from the span of $\{ x_k  \}_{k=1}^n$ to $\ell_1^n$ 
that sends $x_k$ to the  $k^{\text{\, th}}$ unit vector of $\ell_1^n$. 
Set $d_n = \D \{ y_k  \}_{k=1}^n$.  
Reflexivity of $\Y$ gives that $d_n$ tends to infinity.   
 Let $T$ be a (reflexive) Tsirelson-like space 
with normalized unconditional basis $\{ t_n \}$  
such  that for all finite subsets $F$ of natural numbers,  
$$\D \{ t_n \}_{n\in F} \le \max\left\{ 2, \sqrt{d_{\av F \av}}~\right\} \ , 
$$   
where $\av F \av$ is the cardinality of $F$.  
For example, $\{ t_n \}$ can just be an appropriately 
chosen subsequence of the usual basis of the usual Tsirelson space 
[cf\.~CS, Chapter~I].  
Consider the non-completely-continuous map 
$I_{ T } \: T \to \ell_\infty$ that sends  
$t_n$ to the $n^{\text{\, th}}$ 
unit vector of $\ell_\infty$.  By the universality of $I_{\Y}$,  
there 
exists maps $A$ and $B$ such that the following 
diagram commutes.  
$$ 
\vbox{
\settabs\+\indent
&\hskip 1 true in &\hskip 2 true in &\hskip 1 true in \cr 
\+&\hf$T$\hf
  &\hf$\longrightarrow^{^{\hskip-1.25em I_{ T} }}\hskip1em$
 \hf&\hf$\ell_{\infty}$\hf&\cr
\+\cr
\+&\hf$A\uparrow$\hf&  &\hf$\downarrow B$\hf&\cr
\+\cr 
\+&\hf$\Y$\hf
  &\hf$\longrightarrow^{^{\hskip-1.25em I_{\Y}}}\hskip1em$\hf
  &\hf$\ell_\infty$\hf\cr}
$$    

Since   each  $I_{\Y} (y_n)$  is of norm one, 
there exists  $\d >0$ such that 
$\d < \n I_{T } A y_n \n$  for each $n$.  
Each $A y_n$ is of the form 
$$ A y_n = \sum_{m=1}^\infty \a_{n,m}~t_m$$   
and so  there is a sequence $\{ m(n) \}_n$ of natural numbers 
such that $\d < \av \a_{n,m(n)} \av$. Since $\{y_n\}$ tends weakly to
zero, for each $m$ the set of all $n$ for which $m(n)=m$ is finite.
Thus by replacing $\Y$ with the closed span of a suitable subsequence
of $\{y_n\}$, we may assume that the $m(n)$'s are distinct.

Let $T_{*}$ 
be the subspace of $T$ spanned by $\{ t_{m(n)} \}_n$. 
Since $\{y_n\}$ and $\{ t_{m(n)} \} $   
are both unconditional bases, by the diagonalization principle 
[cf\.~LT, Prop\.~1.c.8], the correspondence 
$ y_n \mapsto \a_{n, m(n)} t_{m(n)}$ extends to an
operator  $D\: \Y \to T_{*}$.  
Since $\{ t_{m(n)} \} $ is an unconditional basis 
and $\d < \av \a_{n,m(n)} \av$, the correspondence 
 $ \a_{n, m(n)}~t_{m(n)} \mapsto  t_{m(n)}$ extends to an operator
$M \: T_{*} \to T_{*} $.

By the definition of $d_n$, there exists a sequence $\{ \beta^n_i
\}_{i=1}^n$  such that $\sum_{i=1}^n  \av \beta^n_i \av = 1$ and 
$$ 
\n \sum_{i=1}^n \beta^n_i  y_i \n_{\Y} ~=~ \frac{1}{d_n} \ . 
$$  
   
By the choice of $T$,  for large $n$,  
$$
\frac{1}{\sqrt{d_n}}  ~\le~ 
\n \sum_{i=1}^n \beta^n_i ~ t_{m(i)} \n_{T_{*}} \ . 
$$ 
Since $MD \: \Y \to  T_{*}$ maps   $y_n$ to $t_{m(n)}$, 
$$
\n \sum_{i=1}^n \beta^n_i  t_{m(i)} \n_{T_{*}} 
~\le ~\n MD \n  ~~
\n \sum_{i=1}^n \beta^n_i  y_i \n_{\Y} \ . 
$$ 
This gives that 
$$ 
\frac{1}{\sqrt{d_n}} \le \frac{ \n MD \n }{d_n}  \ , 
$$ 
which cannot be since $d_n$ tends to infinity.   
\hfill\qed 
\enddemo 

The first two paragraphs of the proof of Theorem 5 yield part $(a)$ of
the next proposition.  Part $(b)$ follows from similar considerations
and the Gurarii-James theorem [Ja,~Thm.~2]. 

\proclaim{Proposition~6}   
\roster
\item"{(a)}" Let \ $\fS$ \ be the collection  of all formal identity
operators  into $\ell_\infty$ from  reflexive sequence spaces for which
the unit vectors form a normalized unconditional basis.  Then \ $\fS$ \
is universal for the class of all non-completely-continuous operators.
\item"{(b)}" The collection  
\ $\{ I : \ell_p\to\ell_\infty \, ; \, 1<p<\infty \}$ of formal identity
operators  is universal for the class of all non-completely-continuous
operators whose domain is superreflexive.
\endroster
\endproclaim

Recall that a Banach space $\X$ has the 
Radon-Nikod\'ym Property (RNP) 
[respectively, is strongly regular, has the 
Complete Continuity Property (CCP)]  
if each bounded linear operator from $L_1$ into $\X$ 
is representable [respectively, strongly regular,  completely 
continuous]. The  books [DU], [GGMS], and [T] 
contain  splendid surveys of these properties.  
Here we only  recall that 
a  representable operator is strongly regular and 
a strongly regular operator is completely continuous. 
The first paragraph of the proof of Theorem~1  
uses elementary methods to construct, from an 
operator $T \: L_1 \to \X$ 
that is not completely continuous, a copy 
of $\ell_1$ in the closed span of a norming 
set of $\X$.   
On a much deeper level,  the following fact is well-known.  
\proclaim{Fact}  The following are equivalent.  
\roster 
\item $\ell_1$ embeds into $\X$. 
\item $L_1$ embeds into $\X^*$. 
\item $\X^*$ fails the CCP. 
\item $\X^*$  is not strongly regular. 
\endroster 
\endproclaim 
The  well-known equivalence of (1) and (2) was 
shown by 
Pe\l czy\'nski [P,~for separable~$\X$] and  
Hagler [H, for non-separable~$\X$].   
The other downward implications 
follow from the definitions.  Bourgain 
[B]  used  a 
non-strongly-regular  operator into a dual space  
to construct a copy  of $\ell_1$ in the pre-dual.  
Here the authors wish to  formalize   the following essentially know fact 
which, to the best of our knowledge, has not 
appeared in print as such. 
 
\proclaim{Fact} 
The following are equivalent. 
\roster
\item $\X$ has trivial type. 
\item $\X$ fails super CCP.    
\item $\X$ is not super strongly regular.   
\endroster
\endproclaim  

\demo{Proof}  
To see that (1) implies (2), recall that  
$\X$ has trivial type if and only if 
$\ell_1$ is finitely representable in $\X$ 
and that $L_1$ is finitely representable in $\ell_1$. 
Thus, if $\X$ has trivial type, then   
$L_1$ is finitely representable  in $\X$ 
and so $\X$ cannot  have the super CCP.   
Property (3) formally follows from (2). 
Towards seeing that (3) implies (1), 
consider a space $\X$ that is not 
strongly regular.  From the above fact it follows  that 
$\ell_1$ embeds into $\X^*$.  
Thus $\X^*$ has trivial type, 
which implies  the same for $\X$.  
\hfill\qed   
\enddemo

\widestnumber\no{[DFJP]ZZ}
\def\n #1{\no{[\bf #1]}}

\enddocument